\documentclass[a4paper, 11pt]{article}
\usepackage{amssymb, amsmath,latexsym,amsfonts,amsbsy, amsthm,mathtools,graphicx,CJKutf8,CJKnumb,CJKulem,color,geometry}
\usepackage{verbatim}
\usepackage{hyperref}
\usepackage{newtxtext, newtxmath}
\usepackage[misc]{ifsym}
\usepackage[marginal]{footmisc}
\usepackage[titletoc,title]{appendix}
\def\l{\left}
\def\r{\right}
\def\f{\frac}

\def\az{\alpha}

\def\az{\alpha}

\def\ez{\epsilon}
\def\bz{\beta}
\def\dz{\delta}
\def\gz{\gamma}

\def\beq{\begin{equation}}
\def\eeq{\end{equation}}
\def\be{\begin{equation*}}
\def\ee{\end{equation*}}
\def\beqn{\begin{eqnarray}}
\def\eeqn{\end{eqnarray}}
\def\ben{\begin{eqnarray*}}
\def\een{\end{eqnarray*}}

\theoremstyle{plain}
\theoremstyle{plain}\newtheorem{thm}{Theorem}[section]
\theoremstyle{plain}\newtheorem{prop}{Proposition}[section]
\theoremstyle{plain}\newtheorem{cor}{Corollary}[section]
\theoremstyle{plain}\newtheorem{lem}{Lemma}[section]
\theoremstyle{plain}\newtheorem{rem}{Remark}[section]

\numberwithin{equation}{section}

\makeatletter
\def\thanks#1{\protected@xdef\@thanks{\@thanks
        \protect\footnotetext{#1}}}
\makeatother

\begin{document}

\title{Global solutions of $2$-$D$ cubic Dirac equation with non-compactly supported data}
\author{Qian Zhang}

\date{}

\maketitle

\noindent {\bf{Abstract}}\ \ We are interested in the cubic Dirac equation in two space dimensions. We establish the small data global existence and sharp pointwise decay results for general cubic nonlinearities without additional structure. We also prove the scattering of the Dirac equation for certain classes of nonlinearities. In all the above results we do not require the initial data to have compact support.

\bigskip

\noindent {\bf{Keywords}}\ \ cubic Dirac equation $\cdot$ global-in-time solutions $\cdot$ sharp pointwise decay $\cdot$ ghost weight method

\bigskip

\noindent {\bf{Mathematics Subject Classifications (2010)}}\ \  35L05 $\cdot$ 35L52 $\cdot$ 35L71

\section{Introduction}\label{s1}

Consider the nonlinear Dirac equation in two space dimensions
\beq\label{s1: 1.1}
i\gz^\mu\partial_\mu\psi+m\psi=F(\psi)
\eeq
with initial data
\beq\label{s1: 1.2}
\psi(0,x)=\psi_0(x),
\eeq
where $i\gz^\mu\partial_\mu=i\gz^0\partial_t+i\gz^1\partial_1+i\gz^2\partial_2$ is the Dirac operator, $\partial_a=\partial_{x_a}$ for $a=1,2$, $\psi(t,x):\mathbb{R}^{1+2}\to\mathbb{C}^2$ is a spinor field with mass $m\ge 0$, and $\gz^\mu$ are the Dirac matrices. Dirac matrices are defined by the identities 
\beq\label{s1: 1.3}
\gz^\mu\gz^\nu+\gz^\nu\gz^\mu=-2g^{\mu\nu}I_2,\quad\quad (\gz^\mu)^*=-g_{\mu\nu}\gz^\nu,
\eeq
where $g=\mathrm{diag}(-1,1,1)$ denotes the Minkowski metric in $\mathbb{R}^{1+2}$, $\mu,\nu\in\{0,1,2\}$, $I_2$ is the $2\times 2$ identity matrix and $A^*=(\bar{A})^T$ is the Hermitian conjugate of the matrix $A$. We consider general cubic nonlinearities $F$ and do not require additional structures of $F$, i.e., 
\beq\label{s1: 1.4}
F(\psi)=(\psi^*H\psi)\psi,
\eeq
where $H\in\mathbb{C}^{2\times 2}$ is an arbitrary matrix and $\psi^*$ denotes the complex conjugate transpose of the vector $\psi$. 

In the sequel, we use $C$ to denote a universal constant whose value may change from line to line. As usual, $A\lesssim B$ means that $A\le CB$ for some constant $C$. Given a vector or a scalar $w$ we use Japanese bracket to denote $\langle w\rangle:=(1+|w|^2)^{1/2}$. As usual, we use $\Box=g^{\mu\nu}\partial_\mu\partial_\nu=-\partial_t^2+\partial_1^2+\partial_2^2$ to denote the Minkowski wave operator. 

In quantum field theory the nonlinear Dirac equation is a model of self-interacting Dirac fermions and has been widely used to build relativistic models of extended particles. It was originally formulated in one space dimension known as the Thirring model \cite{T} and in three space dimension known as the Soler model \cite{So}, with cubic nonlinearities $F$ which can be written as 
\beq\label{s1: TS}
F(\psi)=\l\{\begin{array}{l}
(\psi^*\gz^0\gz^\mu\psi)\gz_{\!\mu}\psi,\\
(\psi^*\gz^0\psi)\psi
\end{array}\r.
\eeq
respectively, where $\gz_{\!\mu}=g_{\mu\nu}\gz^\nu$.

In terms of the well-posedness of the Cauchy problem, the scale invariant regularity for the nonlinear Dirac equation in $\mathbb{R}^{1+n}$ is $s_c=\f{n-1}{2}$ and therefore it is expected to be well posed for data $\psi_0\in H^s(\mathbb{R}^n)$ with $s\ge\f{n-1}{2}$. In the low regularity setting, there are numerous results concerning local and global (in time) existence of solutions, see for example \cite{Ev,  Tz, MNN, BH} in the case of three space dimensions. On the other hand, in the case of two space dimensions, Pecher \cite{P} proved the local well-posedness for data in $H^s(\mathbb{R}^2)$ in the almost critical case $s>\f{1}{2}$. Bournaveas and
Candy \cite{BC} proved local well-posedness with initial data in the critical space $H^{\f{1}{2}}(\mathbb{R}^2)$ and global well-posedness for the case $m=0$. Global well-posedness and scattering for the case $m>0$ with small initial data in $H^{\f{1}{2}}(\mathbb{R}^2)$ was established by Bejenaru and Herr \cite{BH2}. We point out that these previous work mostly focused on nonlinearities as in \eqref{s1: TS}.

We study global-in-time existence with pointwise decay of the solution to \eqref{s1: 1.1}-\eqref{s1: 1.2} in $\mathbb{R}^{1+2}$, with nonlinearities as in \eqref{s1: 1.4} and non-compactly supported initial data in weighted Sobolev space of high regularity. We are interested in the massless case $m=0$ since the massive Dirac equation with cubic nonlinearities is easier to treat (see Remark \ref{s1: 1.2}). In a recent work \cite{DL}, the authors considered \eqref{s1: 1.1} with the nonlinearity $F(\psi)=(\psi^*\gz^0\psi)\psi$ and compactly supported initial data, and established global existence and long time dynamics including pointwise decay and scattering, using the hyperboloidal foliation of spacetime. Compared with \cite{DL}, the novelty of our results is that we remove the compactness assumption on the initial data and the structural condition on the nonlinearities. See also \cite{D21} for global existence results on two dimensional coupled wave and Klein-Gordon equations with non-compactly supported initial data.

To conclude, we focus on the the study of the following Cauchy problem in $\mathbb{R}^{1+2}$:
\beq\label{s1: 1.1s}
i\gz^\mu\partial_\mu\psi=F(\psi)=(\psi^*H\psi)\psi,\quad\quad\psi(0,x)=\psi_0(x),
\eeq
for an arbitrary matrix $H\in\mathbb{C}^{2\times 2}$. From now on, we also denote the Dirac operator by
\beq\label{s1: Dirac}
\mathcal{D}:=i\gz^\mu\partial_\mu.
\eeq

\vspace{0.5em}

\noindent${\mathbf{Major\ difficulties\ and\ key\ ideas.}}$ We apply Klainerman's vector field method \cite{K85,K86} to study the Dirac equation \eqref{s1: 1.1s}. Using the identity $\Box=\mathcal{D}^2$, we obtain 
\beq\label{s1: boxd}
\Box \psi=\mathcal{D}^2\psi=\mathcal{D}F(\psi)=\mathcal{D}\l((\psi^*H\psi)\psi\r).
\eeq
We first note that the cubic nonlinearity in \eqref{s1: boxd} violates the standard null condition \cite{K86}. Since the free-linear waves in $\mathbb{R}^{1+2}$ decays at the speed of $\langle t\rangle^{-\f{1}{2}}$,
the identity \eqref{s1: boxd} implies that a linear massless Dirac field in $\mathbb{R}^{1+2}$ should have the same slow decay rate. Hence, the best decay rate we can expect about the $L^2$ norm of the nonlinearity in \eqref{s1: boxd} is 
\be
\|\mathcal{D}\l((\psi^*H\psi)\psi\r)\|_{L^2_x}\lesssim \langle t\rangle^{-1},
\ee
which is the borderline nonintegrable rate. Hence the nonlinearities as in \eqref{s1: 1.4} under consideration may contribute to the long time behavior of the solution to \eqref{s1: 1.1s}. Let us recall the following wave equation in $\mathbb{R}^{1+3}$ with critical nonlinearity
\be
-\Box u=(\partial_tu)^2,\quad\quad (u,\partial_tu)|_{t=0}=(0,u_1)
\ee
for which $u_1$ is compactly supported. John \cite{J81} showed that nontrivial $C^3$ solutions to this equation blow up in finite time. On the other hand, under compactness assumption on the initial data and  additional structure condition on the nonlinearity (i.e. $F(\psi)=(\psi^*\gz^0\psi)\psi$), global existence and unified (in $m$) pointwise decay results for the Dirac equation \eqref{s1: 1.1} were established in \cite{DL}, where the authors use the hyperboloidal foliation method and make full use of the $\gz^0$-structure of $F(\psi)$ to obtain better decay estimates of the solution. For non-compactly supported initial data and nonlinearities as in \eqref{s1: 1.4} without additional structures, our difficulties in using Klainerman's vector field method to study global existence for the problem \eqref{s1: 1.1s} include: $i)$ obtaining good decay estimate for the solution to close the bootstrap argument; $ii)$ dealing with non-compactly supported initial data (in which case the hyperboloidal foliation method cannot be used directly). 

To conquer these difficulties, our key ideas include: $i)$ using the good commutative property of the scaling vector field $L_0$ with the Dirac operator. This means that we can use the full range of the (compatible) vector fields (denoted by $\hat{\Gamma}^I$ for any multi-indices $I$) and obtain the $\langle t-|x|\rangle$ decay of the solution $\psi$ by employing the classical Klainerman-Sobolev inequality; $ii)$ applying Alinhac's ghost weight energy method \cite{Al1} adapted to the Dirac equation and, by a careful calculation when deriving the energy estimate, we obtain the $\gz^0$-structure $F^*\gz^0\psi$ (for nonlinearities $F$ as in \eqref{s1: 1.4}) even if $F$ does not necessarily have this structure itself. This idea of discovering the $\gz^0$-structure is inspired by a recent work \cite{DLMY}, where the authors established global existence, sharp time decay and scattering result for $2D$ Dirac-Klein-Gordon system with non-compactly supported initial data. By a delicate cancellation, we can further write   
\beq\label{s1: Fgzpsi}
F^*\gz^0\psi=[F]_{-}^*\gz^0[\psi]_{+}+[F]_{+}^*\gz^0[\psi]_{-},
\eeq
where $[\psi]_{\pm}=\l(I_2\pm\f{x_a}{|x|}\gz^0\gz^a\r)\psi$ and similarly for $[F]_{\pm}$. The definition \eqref{s1: 1.4} then gives $[F]_{-}=(\psi^*H\psi)[\psi]_{-}$. That is, both terms in \eqref{s1: Fgzpsi} can be written roughly as $[\psi]_{-}\cdot|\psi|^3$. When acting the vector fields $\hat{\Gamma}^I$ on both sides of \eqref{s1: 1.1s} and applying the ghost energy estimate, we obtain the corresponding structure $(\hat{\Gamma}^IF)^*\gz^0\hat{\Gamma}^I\psi$ which can be written roughly as
\beq\label{s1: gzFps}
\l([\hat{\Gamma}^I\psi]_{-}\cdot\psi\cdot\psi+[\psi]_{-}\cdot\hat{\Gamma}^I\psi\cdot\psi\r)\cdot \hat{\Gamma}^I\psi.
\eeq
For the estimate of the first term in \eqref{s1: gzFps}, we need to use the $\langle t-|x|\rangle$ decay of the solution $\psi$, which follows from the Klainerman-Sobolev inequality as stated above. Hence the key to closing the energy estimate is obtaining good pointwise decay estimate of $[\psi]_{-}$ in \eqref{s1: gzFps}. For this, we adopt an idea due to Bournaveas \cite{Bo} and introduce a new function $\Psi$ which solves the wave equation 
$$\Box\Psi=F,\quad\quad (\Psi,\partial_t\Psi)|_{t=0}=(0,-i\gz^0\psi_0),$$ 
and find that we can roughly write $[\psi]_{-}$ as $G_a\Psi$, where $G_a$ denotes the good derivatives. By employing the $L^\infty$ estimate on linear wave equation and using the $\langle t\rangle^{-1}$ decay of  good derivatives, we obtain sufficient decay estimate of $[\psi]_{-}$.

The main result is stated as follows.

\begin{thm}\label{s1: thm1}
Let $N\ge 3$ be an integer. Then there exists $\ez_0>0$ such that for all $0<\ez<\ez_0$ and all initial data $\psi_0$ satisfying the smallness condition
\beq\label{s1: psis}
\sum_{k\le N}\l(\|\langle |x|\rangle^{k+1}\nabla^k\psi_0\|_{L^1_x}+\|\langle |x|\rangle^{k+1}\nabla^k\psi_0\|_{L^2_x}\r)\le\ez,
\eeq
the Cauchy problem \eqref{s1: 1.1s} admits a global-in-time solution $\psi$, which satisfies the following pointwise decay estimate
\ben
|\psi|\lesssim\ez\langle t+|x|\rangle^{-\f{1}{2}}\langle t-|x|\rangle^{-\f{1}{2}}.
\een
\end{thm}

\begin{rem}\label{s1: rem 4}
Theorem \ref{s1: thm1} holds with relaxed condition on the smallness of the initial data (choosing $N$ larger, for example $N\ge 5$), i.e.,
\be
\|\psi_0\|_{L^1_x}+\sum_{k\le N}\|\langle|x|\rangle^{k+1}\nabla^k\psi_0\|_{L^2_x}\le\ez.
\ee
See Appendix \ref{sB} for the proof. 
\end{rem}

\begin{rem}\label{s1: rem3}
Theorem \ref{s1: thm1} also holds for nonlinearities $F=(\psi^*\gz^0\psi)A\psi$, where $A\in\mathbb{C}^{2\times 2}$ is an arbitrary matrix. In this case, the expression of $F$ itself admits a $\gz^0$-structure, and hence $F$ can be written roughly as $[\psi]_{-}|\psi|^2$ as stated in the paragraph above Theorem \ref{s1: thm1}.
\end{rem}

\begin{rem}\label{s1: rem 2}
For the massive Dirac equation
\beq\label{s1: 1.1m}
i\gz^\mu\partial_\mu\psi+m\psi=F(\psi),\quad m>0
\eeq
with initial data \eqref{s1: 1.2} and nonlinearity $F$ as in \eqref{s1: 1.4}, assume without loss of generality that $m=1$. Acting the Dirac operator on both sides of \eqref{s1: 1.1m}, we find that $\psi$ solves the following Klein-Gordon equation
\beq\label{s1: 1.2m}
-\Box \psi+\psi=F(\psi)-i\gz^\mu\partial_\mu F(\psi).
\eeq
Since the nonlinearity $F$ is cubic, using the decay estimates for linear Klein-Gordon equation \cite{Ge} and the bootstrap argument, one can obtain global existence of the solution $\psi$ to \eqref{s1: 1.1m} with small, high-regular, non-compactly supported initial data, and the sharp pointwise decay estimate
\be
|\psi|\lesssim\ez\langle t+|x|\rangle^{-1},
\ee
where $0<\ez\ll 1$ measures the size of the initial data.
\end{rem}

\begin{rem}\label{s1: rem1}
As mentioned above, for compactly supported initial data, the global existence and uniform (in the mass parameter $m\in[0,1]$) pointwise decay estimate for the Cauchy problem in $\mathbb{R}^{1+2}$
\beq\label{s1: 1.1mgz0}
i\gz^\mu\partial_\mu\psi+m\psi=F(\psi)=(\psi^*\gz^0\psi)\psi,\quad\quad\psi(t_0,x)=\psi_0(x),\quad t_0=2
\eeq
were established in \cite{DL}. Our method of treating more general nonlinearities as in \eqref{s1: 1.4} can be adapted to the hyperboloidal foliation case there, and then one can prove the global existence of the solution $\psi$ to the problem
\beq\label{s1: 1.1mH}
i\gz^\mu\partial_\mu\psi+m\psi=F(\psi)=(\psi^*H\psi)\psi,\quad\quad\psi(t_0,x)=\psi_0(x)
\eeq
for an arbitrary matrix $H\in\mathbb{C}^{2\times 2}$ and small compactly supported data $\psi_0$, with the unified pointwise decay estimate 
$$|\psi(t,x)|\lesssim\f{\ez}{t^{\f{1}{2}}(t-|x|)^{\f{1}{2}}+mt}.$$
Namely, the unified (in $m\in[0,1]$) pointwise decay result in \cite[Theorem 1.1]{DL} can be generalized to cubic nonlinearities as  in \eqref{s1: 1.4}. Indeed, the idea used in this paper (see also \cite[Lemma 2.6]{DLMY}) of discovering the $\gz^0$-structure in deriving the ghost energy estimate can be used in proving the hyperboloidal energy estimate (\cite[Proposition 2.1]{DL}) for the Dirac equation (see also \cite[Lemma 2.2]{DLW} and \cite{DW}). 
\end{rem}

Using the ghost weight energy estimate and the pointwise decay result given by Theorem \ref{s1: thm1}, together with the integral formula for the Dirac equation, we obtain the result below concerning the asymptotic behavior (in the Sobolev space) of the global solution obtained in Theorem \ref{s1: thm1}. Precisely, the small global solution scatters as time tends to infinity, and it tends to the solution of a linear Dirac equation in the Sobolev space of high regularity.

\begin{thm}\label{s1: thm2}
Let $N\ge 3, \ez_0>0$ be as in Theorem \ref{s1: thm1} and $\psi_0$ satisfy \eqref{s1: psis} with $0<\ez\le\ez_0$. Suppose $\psi$ is the global solution to \eqref{s1: 1.1s} with $F(\psi)=(\psi^*\gz^0\psi)\psi$. Then the solution $\psi$ scatters linearly as $t\to+\infty$. More precisely, there exists some $\psi^+\in H^N(\mathbb{R}^2)$ such that
\be
\|\psi(t)-S(t)\psi^+\|_{H^N}\le C\langle t\rangle^{-\f{1}{2}}\ln(2+t),\quad\forall t\ge 0
\ee
and
\be
\|\psi(t)-S(t)\psi^+\|_{H^{N-2}}\le C(t)\langle t\rangle^{-\f{1}{2}}\ln(2+t)
\ee
for some $C(t)>0$ satisfying $\lim_{t\to+\infty}C(t)=0$, where $S(t):=e^{-t\gz^0\gz^a\partial_a}$ is the propagator for the linear Dirac equation (see Sect. \ref{s4} for the definition).
\end{thm}

\section{Preliminaries}\label{s2}

\subsection{Notation}\label{s2.1}

We work in the $(1+2)$ dimensional spacetime $\mathbb{R}^{1+2}$ with Minkowski metric $g=(-1,1,1)$, which is used to raise or lower indices. The space indices are denoted by Roman letters $a,b\in\{1,2\}$, and the spacetime indices are denoted by Greek letters $\mu,\nu,\az,\bz\in\{0,1,2\}$. Einstein summation convention for repeated upper and lower indices is adopted throughout the paper. We denote a point in $\mathbb{R}^{1+2}$ by $(t,x)=(x_0,x_1,x_2)$ with $t=x_0,x=(x_1,x_2),x^a=x_a,a=1,2$, and its spacial radius is denoted by $r:=|x|=\sqrt{x_1^2+x_2^2}$. The following vector fields will be used frequently in the analysis:
\begin{itemize}
\item[(i)] Translations: $\partial_\az:=\partial_{x_\az}$, for $\az=0,1,2$.
\item[(ii)] Lorentz boosts: $L_a:=x_a\partial_t+t\partial_a$, for $a=1,2$.
\item[(iii)] Rotation: $\Omega_{12}:=x_1\partial_2-x_2\partial_1$.
\item[(iv)] Scaling: $L_0=t\partial_t+x^a\partial_a$.
\end{itemize}
We also use the modified Lorentz boosts and rotation, first introduced in \cite{Ba},
$$\hat{L}_a:=L_a-\f{1}{2}\gz^0\gz^a,\quad\quad\hat{\Omega}_{12}:=\Omega_{12}-\f{1}{2}\gz^1\gz^2,$$
which enjoy the following commutative property, i.e.
\be
[\hat{L}_a,\mathcal{D}]=[\hat{\Omega}_{12},\mathcal{D}]=0,
\ee
where the commutator $[A,B]$ is defined as 
$$[A,B]:=AB-BA.$$ 
For the ordered sets
$$\{\Gamma_1,\Gamma_2,\cdots,\Gamma_7\}:=\{\partial_0, \partial_1, \partial_2, L_1, L_2, \Omega_{12}, L_0\}$$ 
and
$$\{\hat{\Gamma}_1,\hat{\Gamma}_2,\cdots,\hat{\Gamma}_7\}:=\{\partial_0, \partial_1, \partial_2, \hat{L}_1, \hat{L}_2, \hat{\Omega}_{12}, L_0\}$$ 
and any multi-index $I=(i_1,i_2,\dots,i_7)\in\mathbb{N}^7$ of length $|I|=\sum_{k=1}^7i_k$, we denote
$$\Gamma^I=\prod_{k=1}^7\Gamma_k^{i_k},\quad\mathrm{where}\quad\Gamma=(\Gamma_1,\Gamma_2,\dots,\Gamma_7)$$
and
$$\hat{\Gamma}^I=\prod_{k=1}^7\hat{\Gamma}_k^{i_k},\quad\mathrm{where}\quad\hat{\Gamma}=(\hat{\Gamma}_1,\hat{\Gamma}_2,\dots,\hat{\Gamma}_7).$$
We also introduce the good derivatives
\be
G_a=\f{1}{r}(x_a\partial_t+r\partial_a),\quad\mathrm{for}\;a=1,2.
\ee

\subsection{Estimates on the vector fields and Sobolev inequalities}\label{s2.2}

We first recall the well-known relations
\be
[\Box, \Gamma_k]=0,\quad\mathrm{for}\ k=1,\dots,6,\quad\quad [\Box,L_0]=2\Box.
\ee
Also, by straightforward computation, we have
\beq\label{s2: Dhgz}
[\mathcal{D},\hat{\Gamma}_k]=0,\quad\mathrm{for}\ k=1,\dots,6,\quad\quad [\mathcal{D},L_0]=\mathcal{D}.
\eeq

\begin{lem}\label{s2: 2.2-4}
The following statements hold:
\begin{itemize}
\item[$i)$] 
\beq\label{s2: Dhgz^I}
[\mathcal{D},\hat{\Gamma}^I]=\sum_{|I'|<|I|}c_{I'}\  \mathcal{D}\hat{\Gamma}^{I'}=\sum_{|J'|<|I|}c_{J'}\hat{\Gamma}^{J'}\mathcal{D}
\eeq
for some constants $c_{I'},c_{J'}$.
\item[$ii)$] Let $u=u(t,x):\mathbb{R}^{1+2}\to\mathbb{C}$ be a scalar field and $\Phi=\Phi(t,x):\mathbb{R}^{1+2}\to\mathbb{C}^2$ be a vector field. Then
\beq\label{s2: uPhi}
\hat{\Gamma}^I(u\Phi)=\sum_{I_1+I_2=I}(\Gamma^{I_1}u)(\hat{\Gamma}^{I_2}\Phi).
\eeq
\item[$iii)$] 
\beq\label{s2: hgz^I}
\hat{\Gamma}^I\Phi=\Gamma^I\Phi+\sum_{|I'|<|I|}c_{I'}\Gamma^{I'}\Phi,\quad\quad\Gamma^J\Phi=\hat{\Gamma}^J\Phi+\sum_{|J'|<|J|}c_{J'}\hat{\Gamma}^{J'}\Phi
\eeq
for some constant matrices $c_{I'}, c_{J'}$.
\end{itemize}

\begin{proof}
$i)$ By \eqref{s2: Dhgz}, the equality \eqref{s2: Dhgz^I} holds for $|I|=1$. Assume by induction that \eqref{s2: Dhgz^I} holds for $|I|=l$ with $l\in\mathbb{N}$. For each fixed $|I|=l+1$, we can decompose $I$ into $I_1+I_2$ with $|I_1|=l$ and $|I_2|=1$. Hence,
\ben
[\mathcal{D},\hat{\Gamma}^{I_1+I_2}]&=&[\mathcal{D},\hat{\Gamma}^{I_1}]\hat{\Gamma}^{I_2}+\hat{\Gamma}^{I_1}[\mathcal{D},\hat{\Gamma}^{I_2}]\\
&=&\l(\sum_{|I'_1|<|I_1|}c_{I'_1}\mathcal{D}\hat{\Gamma}^{I'_1}\r)\hat{\Gamma}^{I_2}+c\hat{\Gamma}^{I_1}\mathcal{D}\\
&=&\sum_{|I'_1|<|I_1|}c_{I'_1}\mathcal{D}\hat{\Gamma}^{I'_1}\hat{\Gamma}^{I_2}+c \mathcal{D}\hat{\Gamma}^{I_1}+c\sum_{|J'_1|<|I_1|}c_{J'_1}\mathcal{D}\hat{\Gamma}^{J'_1}\\
&=&\sum_{|J'|<|I|}c_{J'}\mathcal{D}\hat{\Gamma}^{J'}
\een
for some constants $c_{I'_1},c,c_{J'_1},c_{J'}$. This gives the first equality in \eqref{s2: Dhgz^I}. The second equality can be proved similarly.

$ii)$ By straightforward computation,
\be
\hat{L}_a(u\Phi)=(L_au)\Phi+u(\hat{L}_a\Phi)
\ee
and similarly for $\hat{\Omega}_{12}$. Recall from the definitions in Sect. \ref{s2.1} that $\Gamma_k=\hat{\Gamma}_k$ for $k=1,2,3,7$. Hence we have
\be
\hat{\Gamma}_k(u\Phi)=(\Gamma_ku)\Phi+u(\hat{\Gamma}_k\Phi),\quad\mathrm{for}\;k=1,2,\dots,7.
\ee
Assume for some $l\in\mathbb{N}$ and any $|I|\le l$ we have
\be
\hat{\Gamma}^I(u\Phi)=\sum_{I_1+I_2=I}(\Gamma^{I_1}u)(\hat{\Gamma}^{I_2}\Phi).
\ee
For each fixed $|I|=l+1$, we write $I=I_1+I_2$ with $|I_1|=l$ and $|I_2|=1$. Then
\ben
\hat{\Gamma}^{I_2}\hat{\Gamma}^{I_1}(u\Phi)&=&\sum_{J_1+J_2=I_1}\hat{\Gamma}^{I_2}\l((\Gamma^{J_1}u)(\hat{\Gamma}^{J_2}\Phi)\r)\\
&=&\sum_{J_1+J_2=I_1}(\Gamma^{I_2+J_1}u)(\hat{\Gamma}^{J_2}\Phi)+(\Gamma^{J_1}u)(\hat{\Gamma}^{I_2+J_2}\Phi)\\
&=&\sum_{J_1+J_2=I}(\Gamma^{J_1}u)(\hat{\Gamma}^{J_2}\Phi).
\een

$iii)$ Note that \eqref{s2: hgz^I} holds for $|I|=1$. Assume that it holds for some $l\in\mathbb{N}$ and all $|I|\le l$. For each fixed $|I|=l+1$, we write $I=I_1+I_2$ as above, then
\ben
\hat{\Gamma}^{I_2}\hat{\Gamma}^{I_1}\Phi&=&\hat{\Gamma}^{I_2}\Gamma^{I_1}\Phi+\sum_{|I'_1|<|I_1|}c_{I'_1}\hat{\Gamma}^{I_2}\Gamma^{I'_1}\Phi\\
&=&\Gamma^{I_2}\Gamma^{I_1}\Phi+c\Gamma^{I_1}\Phi+\sum_{|I'_1|<|I_1|}c_{I'_1}\l(\Gamma^{I_2}\Gamma^{I'_1}\Phi+\tilde{c}_{I'_1}\Gamma^{I'_1}\Phi\r)\\
&=&\Gamma^I\Phi+\sum_{|I'|<|I|}c_{I'}\Gamma^{I'}\Phi
\een
for some constant matrices $c,c_{I'_1},\tilde{c}_{I'_1},c_{I'}$.
\end{proof}
\end{lem}

We next give the estimate below for the vector fields $\partial_\az$ and good derivatives $G_a$, which will be used in obtaining good pointwise decay result of the solution to \eqref{s1: 1.1s}.

\begin{lem}\label{s2: DMY2.1}
We have
\be
\langle t-r\rangle|\partial u|+\langle t+r\rangle|G_a u|\lesssim\sum_{|I|=1}|\Gamma^Iu|.
\ee
\begin{proof}
The estimate of $|\partial u|$ is well-known, see for example \cite{S}. For $|G_a u|$, we just use the equalities
\be
G_au=\f{1}{r}\l(L_a u+(r-t)\partial_a u\r)=\f{1}{t}\l(L_a u-\f{x_a}{r}(r-t)\partial_t u\r),
\ee
which follow from a direct calculation.
\end{proof}
\end{lem}

Next we present the famous Klainerman-Sobolev inequality whose proof can be found in \cite{A,Ho,S}.

\begin{lem}(See \cite[Proposition 6.5.1]{Ho})\label{s2: K-S}
Let $u=u(t,x)$ be a sufficiently smooth function which decays sufficiently fast at space infinity for each fixed $t\ge 0$. Then for any $t\ge 0$ and $x\in\mathbb{R}^2$, we have 
\be
\langle t+r\rangle^{\f{1}{2}}\langle t-r\rangle^{\f{1}{2}}|u(t,x)|\lesssim\sum_{|I|\le 2}\|\Gamma^Iu(t,\cdot)\|_{L^2_x}.
\ee
\end{lem}

\subsection{Reformulation of the solution to \eqref{s1: 1.1s}}\label{s2.3}

Below we reformulate the solution to the massless Dirac equation $\psi$ as the Dirac operator acting on the solution $\Psi$ to a linear wave equation. Then by decomposing $\psi$ into two parts $[\psi]_{-}$ and $[\psi]_{+}$ and writing $[\psi]_{-}$ as the good derivatives $G_a\Psi$, we obtain good pointwise decay estimate for $[\psi]_{-}$. This is an important observation and will be used in closing the bootstrap estimate for the ghost weight energy (see Sect. \ref{ss3.1}).

For any $\mathbb{C}^2$-valued function $\Phi$, we define the Hermitian matrices
\be
T_{-}:=I_2-\f{x_a}{r}\gz^0\gz^a,\quad\quad T_{+}:=I_2+\f{x_a}{r}\gz^0\gz^a,
\ee
and let
\beq\label{s2.3: Phi-}
[\Phi]_{-}=T_{-}\Phi=\Phi-\f{x_a}{r}\gz^0\gz^a\Phi,\quad\quad [\Phi]_{+}=T_{+}\Phi=\Phi+\f{x_a}{r}\gz^0\gz^a\Phi.
\eeq

\begin{lem}\label{s2: DL4.2}
The following statements hold:
\begin{itemize}
\item[$i)$] Let $\varphi,\Phi$ be two $\mathbb{C}^2$-valued functions, then
\beq\label{s2: Psi*Phi}
\varphi^{*}\gz^0\Phi=\f{[\varphi]^*_{-}\gz^0[\Phi]_{+}+[\varphi]^*_{+}\gz^0[\Phi]_{-}}{4}.
\eeq
\item[$ii)$] For any $\mathbb{C}^2$-valued function $\Phi$, let $\varphi:=\mathcal{D}\Phi$. Then 
\beq\label{s2: Phi-}
[\varphi]_{-}=i\l(I_2-\f{x_b}{r}\gz^0\gz^b\r)\gz^\mu\partial_\mu\Phi=i\l(I_2-\f{x_b}{r}\gz^0\gz^b\r)\gz^aG_a\Phi.
\eeq

\end{itemize}

\begin{proof}
$i)$ We write
\be
\varphi=\f{[\varphi]_{-}+[\varphi]_{+}}{2},\quad\quad \Phi=\f{[\Phi]_{-}+[\Phi]_{+}}{2}.
\ee
By direct computation,
\beqn\label{s2: T-}
T_{-}\gz^0T_{-}&=&\l(I_2-\f{x_b}{r}\gz^0\gz^b\r)\gz^0\l(I_2-\f{x_a}{r}\gz^0\gz^a\r)\nonumber\\
&=&\gz^0-\f{x_b}{r}\gz^0\gz^b\gz^0-\f{x_a}{r}\gz^0\gz^0\gz^a+\f{x_ax_b}{r^2}\gz^0\gz^b\gz^0\gz^0\gz^a\nonumber\\
&=&\gz^0+\f{x_a^2}{r^2}\gz^0(\gz^a)^2+\f{x_1x_2}{r^2}\gz^0\l(\gz^1\gz^2+\gz^2\gz^1\r)\nonumber\\
&=&0,
\eeqn
and similarly,
\beqn\label{s2: T+}
T_{+}\gz^0T_{+}&=&\l(I_2+\f{x_b}{r}\gz^0\gz^b\r)\gz^0\l(I_2+\f{x_a}{r}\gz^0\gz^a\r)\nonumber\\
&=&\gz^0+\f{x_b}{r}\gz^0\gz^b\gz^0+\f{x_a}{r}\gz^0\gz^0\gz^a+\f{x_ax_b}{r^2}\gz^0\gz^b\gz^0\gz^0\gz^a\nonumber\\
&=&0.
\eeqn
It follows that
\be
[\varphi]_{-}^*\gz^0[\Phi]_{-}=\varphi^*T_{-}\gz^0T_{-}\Phi=0,\quad\quad [\varphi]_{+}^*\gz^0[\Phi]_{+}=\varphi^*T_{+}\gz^0T_{+}\Phi=0,
\ee
which implies $i)$. 

$ii)$ Using the relation
\be
\partial_a=G_a-\f{x_a}{r}\partial_t,
\ee
we can write
\be
\gz^0\partial_t+\gz^a\partial_a=\gz^0\l(I_2-\f{x_a}{r}\gz^0\gz^a\r)\partial_t+\gz^aG_a.
\ee
By \eqref{s2: T-},
\be
\l(I_2-\f{x_b}{r}\gz^0\gz^b\r)\gz^0\l(I_2-\f{x_a}{r}\gz^0\gz^a\r)=T_{-}\gz^0T_{-}=0,
\ee
hence \eqref{s2: Phi-} holds.

\end{proof}
\end{lem}

Let $\psi$ be the solution to \eqref{s1: 1.1s}. We choose $\Psi$ which solves the problem
\beq\label{s3: DL5.3}
\Box\Psi=\mathcal{D}\psi=F(\psi),\quad\quad \l(\Psi,\partial_t\Psi\r)|_{t=0}=\l(0,-i\gz^0\psi_0\r),
\eeq
where $F$ is as in \eqref{s1: 1.1s}. Then 
\beq\label{s2: psi=dPsi}
\psi=\mathcal{D}\Psi=i\gz^\mu\partial_\mu\Psi,
\eeq
since $\varphi:=\psi-\mathcal{D}\Psi$ verifies 
$$\mathcal{D}\varphi=0,\quad\quad\varphi(0,x)=0.$$ 
Using \eqref{s2: Phi-} we can write
\beq\label{s2: psi-} 
[\psi]_{-}=i\l(I_2-\f{x_b}{r}\gz^0\gz^b\r)\gz^aG_a\Psi.
\eeq
In addition, we have the following estimate on $|[\hat{\Gamma}^I\psi]_{-}|$.

\begin{lem}\label{s2: lem psi-}
Let $\psi$ be the solution to \eqref{s1: 1.1s} and $\Psi$ be chosen to satisfy \eqref{s3: DL5.3}. Then we have
\be
|[\hat{\Gamma}^I\psi]_{-}|\lesssim\langle t+r\rangle^{-1}\sum_{|I'|\le|I|+1}|\Gamma^{I'}\Psi|\lesssim\langle t-r\rangle^{-\f{1}{2}}\langle t+r\rangle^{-\f{3}{2}}\sum_{|I'|\le |I|+3}\|\Gamma^{I'}\Psi\|_{L^2_x}
\ee
and
\be
|\hat{\Gamma}^I\psi|\lesssim\langle t-r\rangle^{-1}\sum_{|I'|\le|I|+1}|\Gamma^{I'}\Psi|\lesssim\langle t-r\rangle^{-\f{3}{2}}\langle t+r\rangle^{-\f{1}{2}}\sum_{|I'|\le |I|+3}\|\Gamma^{I'}\Psi\|_{L^2_x}.
\ee
\begin{proof}
Acting the vector field $\hat{\Gamma}^I$ on both sides of \eqref{s2: psi=dPsi} and using \eqref{s2: Dhgz^I}, we obtain
\be
\hat{\Gamma}^I\psi=\mathcal{D}\hat{\Gamma}^I\Psi+\sum_{|I'|<|I|}c_{I'}\ \mathcal{D}\hat{\Gamma}^{I'}\Psi
\ee
for some constants $c_{I'}$. Hence \eqref{s2: Phi-} implies
\beq\label{s2: hgz^Ipsi-} 
[\hat{\Gamma}^I\psi]_{-}=i\l(I_2-\f{x_b}{r}\gz^0\gz^b\r)\l(\gz^aG_a\hat{\Gamma}^I\Psi+\sum_{|I'|<|I|}c_{I'}\gz^aG_a\hat{\Gamma}^{I'}\Psi\r).
\eeq
By Lemma \ref{s2: DMY2.1} and \eqref{s2: hgz^Ipsi-}, we have
\be
|[\hat{\Gamma}^I\psi]_{-}|\lesssim\sum_{|I'|\le|I|}|G_a\hat{\Gamma}^{I'}\Psi|\lesssim\langle t+r\rangle^{-1}\sum_{|I'|\le|I|,|J|=1}|\Gamma^J\hat{\Gamma}^{I'}\Psi|\lesssim\langle t+r\rangle^{-1}\sum_{|K|\le|I|+1}|\Gamma^K\Psi|
\ee
and 
\be
|\hat{\Gamma}^I\psi|\lesssim\sum_{|I'|\le|I|}|\partial\hat{\Gamma}^{I'}\Psi|\lesssim\langle t-r\rangle^{-1}\sum_{|I'|\le|I|,|J|=1}|\Gamma^J\hat{\Gamma}^{I'}\Psi|\lesssim\langle t-r\rangle^{-1}\sum_{|K|\le|I|+1}|\Gamma^K\Psi|,
\ee
where we use \eqref{s2: hgz^I} in these inequalities. By Lemma \ref{s2: K-S},
\be
|\Gamma^K\Psi|\lesssim\langle t-r\rangle^{-\f{1}{2}}\langle t+r\rangle^{-\f{1}{2}}\sum_{|J|\le 2}\|\Gamma^J\Gamma^K\Psi\|_{L^2_x}\lesssim\langle t-r\rangle^{-\f{1}{2}}\langle t+r\rangle^{-\f{1}{2}}\sum_{|K'|\le |K|+2}\|\Gamma^{K'}\Psi\|_{L^2_x}.
\ee
The conclusion follows from the last three estimates.
\end{proof}
\end{lem}

Recall that $\Psi$ is the solution to \eqref{s3: DL5.3}. The lemma below gives $L^\infty$ estimate on $2D$ linear wave equation.

\begin{lem}(See \cite[Theorems 4.6.1, 4.6.2]{LZ})\label{s2: DMY2.78}
Let $u$ be the solution to the Cauchy problem in $\mathbb{R}^{1+2}$
\be
\l\{\begin{array}{rcl}
-\Box u(t,x)&=&f(t,x),\\
(u,\partial_tu)|_{t=0}&=&(u_0,u_1).
\end{array}\r.
\ee
Then we have  
\be
\|u(t,x)\|_{L^\infty_x}\lesssim\langle t\rangle^{-\f{1}{2}}\Bigg\{\|u_0\|_{W^{2,1}}+\|u_1\|_{W^{1,1}}+
\sum_{|I|\le 1}\int_0^t(1+\tau)^{-\f{1}{2}}\|\Gamma^If(\tau,x)\|_{L^1_x}\rm{d}\tau\Bigg\}.
\ee

\end{lem}

\section{Proof of Theorem \ref{s1: thm1}}\label{s3}

\subsection{Ghost weight energy estimate}\label{ss3.1}

To treat general nonlinearities as in \eqref{s1: 1.4} without additional structure, we prove the ghost weight energy estimate below, which generates the $\gz^0$ structure $F^*\gz^0\psi$ on the right hand side. This brings great advantages since we can decompose $F^*\gz^0\psi$ into terms involving $[\psi]_{-}$ (using Lemma \ref{s2: DL4.2}) and therefore applying the decay estimates given by Lemmas \ref{s2: lem psi-} and \ref{s2: DMY2.78} to close the ghost energy estimate in the bootstrap assumption.

\begin{prop}\label{s3: Ghost}
Let $\psi$ be the solution to the Dirac equation \eqref{s1: 1.1s}, then for any $t\ge 0$, we have the following ghost weight energy estimate:
\be
E^D_{gst}(t,\psi)\lesssim \|\psi(0)\|_{L^2_x}^2+\int_0^t\|F^*\gz^0\psi\|_{L^1_x}\rm{d}\tau,
\ee
where 
\beq\label{s3: G1''}
E^D_{gst}(t,\psi):=\|\psi(t)\|_{L^2_x}^2+\int_0^t\l\|\f{[\psi]_-}{\langle r-\tau\rangle^{\f{1+\dz}{2}}}\r\|_{L^2_x}^2\rm{d}\tau.
\eeq
Here $[\psi]_-=\psi-\f{x_a}{r}\gz^0\gz^a\psi$ is defined as in \eqref{s2.3: Phi-} and $\dz>0$ is a constant.
\begin{proof}
Let $q(t,x)=\tilde{q}(r-t)$, where
\be
\tilde{q}(s)=\int_{-\infty}^s\f{1}{\langle\tau\rangle^{1+\dz}}\rm{d}\tau.
\ee
Multiplying on both sides of \eqref{s1: 1.1} by $-ie^q\psi^*\gz^0$, we obtain
\be
e^q\psi^*\partial_t\psi+e^q\psi^*\gz^0\gz^a\partial_a\psi=-ie^q\psi^*\gz^0F.
\ee
Taking the complex conjugate of the last equality, we find that
\be
e^q\partial_t\psi^*\psi+e^q\partial_a\psi^*\gz^0\gz^a\psi=ie^qF^*\gz^0\psi.
\ee
Summarizing the above two equalities and using Leibniz rule, we derive
\beq\label{s3: DL4.3}
\partial_t\l(e^q\psi^*\psi\r)+\partial_a\l(e^q\psi^*\gz^0\gz^a\psi\r)-(\partial_tq) e^q\psi^*\psi-(\partial_aq) e^q\psi^*\gz^0\gz^a\psi=G,
\eeq
where 
\be
G:=ie^q\l(F^*\gz^0\psi-\psi^*\gz^0F\r).
\ee
Substituting the relations
\be
\partial_tq=-\tilde{q}'(r-t)=-\f{1}{\langle r-t\rangle^{1+\dz}},\quad\quad\partial_aq=\tilde{q}'(r-t)=\f{1}{\langle r-t\rangle^{1+\dz}}\f{x_a}{r}
\ee
into \eqref{s3: DL4.3}, we get
\beq\label{s3: DL4.3'}
\partial_t\l(e^q\psi^*\psi\r)+\partial_a\l(e^q\psi^*\gz^0\gz^a\psi\r)+\f{e^q}{\langle r-t\rangle^{1+\dz}}\l(\psi^*\psi-\f{x_a}{r}\psi^*\gz^0\gz^a\psi\r)=G.
\eeq
By straightforward computation,
\be
\l(\psi-\f{x_a}{r}\gz^0\gz^a\psi\r)^*\l(\psi-\f{x_a}{r}\gz^0\gz^a\psi\r)=2\l(\psi^*\psi-\f{x_a}{r}\psi^*\gz^0\gz^a\psi\r).
\ee
Hence \eqref{s3: DL4.3'} yields
\beq\label{s3: DL4.3''}
\partial_t\l(e^q\psi^*\psi\r)+\partial_a\l(e^q\psi^*\gz^0\gz^a\psi\r)+\f{e^q}{2\langle r-t\rangle^{1+\dz}}\l(\psi-\f{x_a}{r}\gz^0\gz^a\psi\r)^*\l(\psi-\f{x_a}{r}\gz^0\gz^a\psi\r)=G.
\eeq
Integrating \eqref{s3: DL4.3''} over $\mathbb{R}^2$, we obtain
\beq\label{s3: G1}
\partial_t\|e^{q/2}\psi\|_{L^2_x}^2+\f{1}{2}\l\|e^{q/2}\f{\psi-\f{x_a}{r}\gz^0\gz^a\psi}{\langle r-t\rangle^{\f{1+\dz}{2}}}\r\|_{L^2_x}^2\le 2\|e^qF^*\gz^0\psi\|_{L^1_x}.
\eeq
Integrating \eqref{s3: G1} over $[0,t]$, we further obtain
\beq\label{s3: G1'}
E^D_{gst}(t,\psi)\lesssim \|\psi(0)\|_{L^2_x}^2+\int_0^t\|F^*\gz^0\psi\|_{L^1_x}\rm{d}\tau,
\eeq
where we use that $e^q\sim 1$.
\end{proof}
\end{prop}

\subsection{Proof of Theorem \ref{s1: thm1}}\label{ss3.2}

\noindent{$\mathbf{Bootstrap\ setting.}$} Let the assumptions in Theorem \ref{s1: thm1} hold. Fix $N\ge 3$ and $0<\dz\ll 1$. Following the local well-posedness theory for the Dirac equation, there exist constants $C_0>0$ and $T>0$ (small) such that \eqref{s1: 1.1s} admits a solution in $[0,T)$ with 
\be
\sum_{|I|\le N}[E^D_{gst}(0,\hat{\Gamma}^I\psi)]^{\f{1}{2}}+\sum_{|I|\le N-2}\|\hat{\Gamma}^I\psi(0,x)\|_{L^\infty_x}\le C_0\ez.
\ee
In addition, let $\Psi$ be the solution to \eqref{s3: DL5.3}, i.e.,
\beq\label{s3: DL5.3'}
\Box\Psi=F(\psi)=(\psi^*H\psi)\psi,\quad\quad \l(\Psi,\partial_t\Psi\r)|_{t=0}=\l(0,-i\gz^0\psi_0\r),
\eeq
where $H$ is as in \eqref{s1: 1.1s}, then we also have
\be
\sum_{|I|\le N-1}\l(\|\Gamma^I\Psi(0)\|_{W^{2,1}}+\|\partial_t\Gamma^I\Psi(0)\|_{W^{1,1}}\r)\le C_0\ez.
\ee
Let $C_1\gg 1$ and $0<\ez_0\ll C_1^{-1}$ be chosen later. We assume the following bootstrap setting:
\beq\label{s3: bs}
\sum_{|I|\le N}[E^D_{gst}(t,\hat{\Gamma}^I\psi)]^{\f{1}{2}}+\sum_{|I|\le N-2}\langle t\rangle^{\f{3}{2}-\dz}|[\hat{\Gamma}^I\psi]_{-}|\le C_1\ez,
\eeq
where $\ez\le\ez_0$ measures the size of the initial data, and (see \eqref{s3: G1''})
\beq\label{s3: ghost}
E^D_{gst}(t,\hat{\Gamma}^I\psi):=\|\hat{\Gamma}^I\psi(t)\|_{L^2_x}^2+\int_0^t\l\|\f{[\hat{\Gamma}^I\psi]_{-}}{\langle r-\tau\rangle^{\f{1+\dz}{2}}}\r\|_{L^2_x}^2\rm{d}\tau.
\eeq
Define
\beq\label{s3: T_*}
T_*=\sup\{t\in(0,\infty): \psi\  \mathrm{satisfies\ }\eqref{s3: bs}\ \mathrm{on}\ [0,t]\}.
\eeq

Theorem \ref{s1: thm1} follows from the result below.

\begin{prop}\label{s3: prop 1}
For all initial data $\psi_0$ satisfying the assumption in Theorem \ref{s1: thm1}, we have $T_*=\infty$.
\end{prop}

Below we give the proof of Proposition \ref{s3: prop 1}. We emphasize that the implied constants in $\lesssim$ do not depend on the constants $C_1$ and $\ez$ appearing in the bootstrap assumption \eqref{s3: bs}. 

We first observe that, by the bootstrap setting \eqref{s3: bs} and Lemma \ref{s2: K-S}, the estimate
\beq\label{s3: gzpsiinf}
\langle t+r\rangle^{\f{1}{2}}\langle t-r\rangle^{\f{1}{2}}\sum_{|I'|\le N-2}|\Gamma^{I'}\psi(t,x)|\lesssim\sum_{|I|\le N}\|\Gamma^I\psi(t)\|_{L^2_x}\lesssim C_1\ez
\eeq
holds true, where we use \eqref{s2: hgz^I} in the last inequality. \\

\noindent$\textit{Proof\ of\ Proposition\ \ref{s3: prop 1}:}$ We divide the proof into two steps.\\

\noindent{$\mathbf{Step\ 1.}$} Refining the estimate of $E^D_{gst}(t,\hat{\Gamma}^I\psi)$. Acting the vector field $\hat{\Gamma}^I$ on both sides of \eqref{s1: 1.1s} and using \eqref{s2: Dhgz^I}, we obtain
\beq\label{s3: Dhgzpsi}
\mathcal{D}\hat{\Gamma}^I\psi=\hat{\Gamma}^IF+\sum_{|I'|<|I|}c_{I'}\hat{\Gamma}^{I'}F
\eeq
for some constants $c_{I'}$. Applying Proposition \ref{s3: Ghost} to $\hat{\Gamma}^I\psi, |I|\le N$ and using \eqref{s3: Dhgzpsi}, we obtain
\be
E^D_{gst}(t,\hat{\Gamma}^I\psi)\lesssim \|\hat{\Gamma}^I\psi(0)\|_{L^2_x}^2+\sum_{|I'|\le|I|}\int_0^t\|(\hat{\Gamma}^{I'}F)^*\gz^0\hat{\Gamma}^I\psi\|_{L^1_x}\rm{d}\tau.
\ee
For each $|I'|\le|I|$, using Lemma \ref{s2: DL4.2}, we have
\beqn\label{s3: Fgz0psi}
|(\hat{\Gamma}^{I'}F)^*\gz^0\hat{\Gamma}^I\psi|&\lesssim&\l|[\hat{\Gamma}^{I'}F]^*_{-}\gz^0[\hat{\Gamma}^I\psi]_{+}+[\hat{\Gamma}^{I'}F]^*_{+}\gz^0[\hat{\Gamma}^I\psi]_{-}\r|\nonumber\\
&\lesssim&\l|\langle r-\tau\rangle^{\f{1+\dz}{2}}\hat{\Gamma}^{I'}F\r|\cdot\l|\f{[\hat{\Gamma}^I\psi]_{-}}{\langle r-\tau\rangle^{\f{1+\dz}{2}}}\r|+|[\hat{\Gamma}^{I'}F]_{-}|\cdot|\hat{\Gamma}^I\psi|.
\eeqn
By \eqref{s2: uPhi}, we have
\beq\label{s3: hgzF}
\hat{\Gamma}^{I'}F=\hat{\Gamma}^{I'}\l((\psi^*H\psi)\psi\r)=\sum_{I'_1+I'_2=I'}(\Gamma^{I'_1}(\psi^*H\psi))(\hat{\Gamma}^{I'_2}\psi)=\sum_{I'_1+I'_2+I'_3=I'}(\Gamma^{I'_1}\psi)^*H(\Gamma^{I'_2}\psi)(\hat{\Gamma}^{I'_3}\psi),
\eeq
which implies
\beq\label{s3: hgzF-}
[\hat{\Gamma}^{I'}F]_{-}=\sum_{I'_1+I'_2+I'_3=I'}(\Gamma^{I'_1}\psi)^*H(\Gamma^{I'_2}\psi)[\hat{\Gamma}^{I'_3}\psi]_{-}.
\eeq
Substituting \eqref{s3: hgzF} and \eqref{s3: hgzF-} into \eqref{s3: Fgz0psi}, using \eqref{s2: hgz^I}, and recalling that $|I|\le N$, we get
\beqn\label{s3: Fgz0psi'}
|(\hat{\Gamma}^{I'}F)^*\gz^0\hat{\Gamma}^I\psi|&\lesssim&\sum_{|I_1|+|I_2|+|I_3|\le|I|}\langle r-\tau\rangle^{\f{1+\dz}{2}}|\Gamma^{I_1}\psi|\cdot|\Gamma^{I_2}\psi|\cdot|\Gamma^{I_3}\psi|\cdot\l|\f{[\hat{\Gamma}^I\psi]_{-}}{\langle r-\tau\rangle^{\f{1+\dz}{2}}}\r|\nonumber\\
&+&\sum_{\substack{|I_1|+|I_2|+|I_3|\le|I|,\\|J|\le |I|}}|\Gamma^{I_1}\psi|\cdot|\Gamma^{I_2}\psi|\cdot|[\hat{\Gamma}^{I_3}\psi]_{-}|\cdot|\Gamma^J\psi|\nonumber\\
&\lesssim&\sum_{\substack{|I_1|,|I_2|\le N-2,\\|I_3|\le N}}\|\langle r-\tau\rangle^{\f{1+\dz}{2}}\Gamma^{I_1}\psi\|_{L^\infty_x}\cdot\|\Gamma^{I_2}\psi\|_{L^\infty_x}\cdot|\Gamma^{I_3}\psi|\cdot\l|\f{[\hat{\Gamma}^I\psi]_{-}}{\langle r-\tau\rangle^{\f{1+\dz}{2}}}\r|\nonumber\\
&+&\sum_{\substack{|I_1|,|I_2|\le N-2,\\|I_3|,|J|\le N}}\|\langle r-\tau\rangle^{\f{1+\dz}{2}}\Gamma^{I_1}\psi\|_{L^\infty_x}\cdot\|\Gamma^{I_2}\psi\|_{L^\infty_x}\cdot\l|\f{[\hat{\Gamma}^{I_3}\psi]_{-}}{\langle r-\tau\rangle^{\f{1+\dz}{2}}}\r|\cdot|\Gamma^J\psi|\nonumber\\
&+&\sum_{\substack{|I_2|,|I_3|\le N-2,\\|I_1|,|J|\le N}}\|\Gamma^{I_2}\psi\|_{L^\infty_x}\cdot\|[\hat{\Gamma}^{I_3}\psi]_{-}\|_{L^\infty_x}\cdot|\Gamma^{I_1}\psi|\cdot|\Gamma^J\psi|\nonumber\\
&\lesssim&(C_1\ez)^2\sum_{|I_1|,|I_2|\le N}\l(\langle \tau\rangle^{-1+\f{\dz}{2}}|\Gamma^{I_1}\psi|\cdot\l|\f{[\hat{\Gamma}^{I_2}\psi]_{-}}{\langle r-\tau\rangle^{\f{1+\dz}{2}}}\r|+\langle \tau\rangle^{-2+\dz}|\Gamma^{I_1}\psi|\cdot|\Gamma^{I_2}\psi|\r),\nonumber\\
\eeqn
where we use \eqref{s3: gzpsiinf}, the bootstrap setting \eqref{s3: bs}, and recall that $N\ge 3$. It follows that 
\ben
&&\|(\hat{\Gamma}^{I'}F)^*\gz^0\hat{\Gamma}^I\psi\|_{L^1_x}\\
&\lesssim&(C_1\ez)^2\sum_{|I_1|,|I_2|\le N}\l(\langle \tau\rangle^{-1+\f{\dz}{2}}\|\Gamma^{I_1}\psi\|_{L^2_x}\cdot\l\|\f{[\hat{\Gamma}^{I_2}\psi]_{-}}{\langle r-\tau\rangle^{\f{1+\dz}{2}}}\r\|_{L^2_x}+\langle \tau\rangle^{-2+\dz}\|\Gamma^{I_1}\psi\|_{L^2_x}\cdot\|\Gamma^{I_2}\psi\|_{L^2_x}\r)\\
&\lesssim&(C_1\ez)^3\sum_{|J|\le N}\l(\langle \tau\rangle^{-1+\f{\dz}{2}}\l\|\f{[\hat{\Gamma}^{J}\psi]_{-}}{\langle r-\tau\rangle^{\f{1+\dz}{2}}}\r\|_{L^2_x}+C_1\ez\langle \tau\rangle^{-2+\dz}\r).
\een
where we use\eqref{s2: hgz^I} and the bootstrap assumption \eqref{s3: bs}. Using \eqref{s3: bs} again, we obtain
\ben
&&\sum_{|I'|\le|I|}\int_0^t\|(\hat{\Gamma}^{I'}F)^*\gz^0\hat{\Gamma}^I\psi\|_{L^1_x}\rm{d}\tau\\
&\lesssim&(C_1\ez)^3\sum_{|J|\le N}\l[\l(\int_0^t\langle \tau\rangle^{-2+\dz}\rm{d}\tau\r)^{\f{1}{2}}\l(\int_0^t\l\|\f{[\hat{\Gamma}^{J}\psi]_{-}}{\langle r-\tau\rangle^{\f{1+\dz}{2}}}\r\|^2_{L^2_x}\rm{d}\tau\r)^{\f{1}{2}}+C_1\ez\int_0^t\langle \tau\rangle^{-2+\dz}\rm{d}\tau\r]\\
&\lesssim&(C_1\ez)^4.
\een
We conclude that
\be
\sum_{|I|\le N}[E^D_{gst}(t,\hat{\Gamma}^I\psi)]^{\f{1}{2}}\le C_0\ez+C(C_1\ez)^2.
\ee
This strictly improves the bootstrap estimate of $E^D_{gst}(t,\hat{\Gamma}^I\psi)$ in \eqref{s3: bs} for $C_1$ sufficiently large and $\ez$ sufficiently small (choose $C_1\ge 8C_0$ and then $C(C_1\ez)\le 1/8$).

\vspace{0.5em}

\noindent{$\mathbf{Step\ 2.}$} Refining the estimate of $[\hat{\Gamma}^I\psi]_{-}$. Let $\Psi$ be chosen to satisfy \eqref{s3: DL5.3'}. Acting the vector fields $\Gamma^I, |I|\le N-1$ on both sides of \eqref{s3: DL5.3'}, we derive
\beq\label{s3: bgzIPsi}
\Box \Gamma^I\Psi=\Gamma^IF+\sum_{|J'|<|I|}c_{J'}\Gamma^{J'}F
\eeq
for some constants $c_{J'}$. By Lemma \ref{s2: DMY2.78}, we have
\beqn\label{s3: gPsiB}
\|\Gamma^I\Psi(t)\|_{L^\infty_x}&\lesssim&\langle t\rangle^{-\f{1}{2}}\l(\|\Gamma^I\Psi(0)\|_{W^{2,1}}+\|\partial_t\Gamma^I\Psi(0)\|_{W^{1,1}}+\sum_{|I'|\le|I|+1}\int_0^t\langle\tau\rangle^{-\f{1}{2}}\|\Gamma^{I'}F(\tau)\|_{L^1_x}\rm{d}\tau\r).
\eeqn
Since
\be
\Gamma^{I'}F=\Gamma^{I'}\l((\psi^*H\psi)\psi\r)=\sum_{I'_1+I'_2=I'}(\Gamma^{I'_1}(\psi^*H\psi))(\Gamma^{I'_2}\psi)=\sum_{I'_1+I'_2+I'_3=I'}(\Gamma^{I'_1}\psi)^*H(\Gamma^{I'_2}\psi)(\Gamma^{I'_3}\psi),
\ee
we have
\ben
\|\Gamma^{I'}F(\tau)\|_{L^1_x}&\lesssim&\sum_{|I_1|+|I_2|+|I_3|\le N}\||\Gamma^{I_1}\psi|\cdot|\Gamma^{I_2}\psi|\cdot|\Gamma^{I_3}\psi|\|_{L^1_x}\\
&\lesssim&\sum_{\substack{|I_1|\le N-2,\\|I_2|,|I_3|\le N}}\|\Gamma^{I_1}\psi\|_{L^\infty_x}\|\Gamma^{I_2}\psi\|_{L^2_x}\|\Gamma^{I_3}\psi\|_{L^2_x}\lesssim (C_1\ez)^3\langle\tau\rangle^{-\f{1}{2}},
\een
where we use \eqref{s3: gzpsiinf} and \eqref{s3: bs} in the last inequality. If follows that
\be
\|\Gamma^I\Psi(t)\|_{L^\infty_x}\le\l(C_0\ez+C(C_1\ez)^3\ln(2+t)\r)\langle t\rangle^{-\f{1}{2}},\quad\mathrm{for\;any\;}|I|\le N-1.
\ee
By Lemma \ref{s2: lem psi-}, we have
\beq\label{s3.2: gpsi-}
\sum_{|I'|\le N-2}|[\hat{\Gamma}^{I'}\psi]_{-}|\lesssim\langle t\rangle^{-1}\sum_{|I|\le N-1}\|\Gamma^{I}\Psi\|_{L^\infty_x}
\le C\l(C_0\ez+C(C_1\ez)^3\r)\langle t\rangle^{-\f{3}{2}}\ln(2+t).
\eeq
Note that we also obtain from Lemma \ref{s2: lem psi-} that
\beq\label{s3.2: gpsi}
\sum_{|I'|\le N-2}|\hat{\Gamma}^{I'}\psi|\lesssim\langle t-r\rangle^{-1}\sum_{|I|\le N-1}\|\Gamma^{I}\Psi\|_{L^\infty_x}
\lesssim\l(\ez+(C_1\ez)^3\r)\langle t-r\rangle^{-1}\langle t\rangle^{-\f{1}{2}}\ln(2+t).
\eeq

In particular, the estimate \eqref{s3.2: gpsi-} strictly improves the bootstrap estimate of $\sum_{|I|\le N-2}|[\hat{\Gamma}^I\psi]_{-}|$ in \eqref{s3: bs} for $C_1$ sufficiently large and $\ez$ sufficiently small (choose $C_1\ge 8CC_0$ and then $C^2(C_1\ez)^2\le 1/8$).

In conclusion, for all initial data $\psi_0$ satisfying the assumption in Theorem \ref{s1: thm1}, we show that $T_*=\infty$ and hence the proof of Proposition \ref{s3: prop 1} is complete.

\section{Scattering for the Dirac field}\label{s4}

In this section we briefly discuss about the scattering of the Dirac equation \eqref{s1: 1.1s}. We show that for the nonlinearity $F=(\psi^*\gz^0\psi)\psi$, the Dirac field $\psi$ scatters linearly in the Sobolev space of high regularity.

We need the following result from \cite[Theorem 4.4]{DL}, which provides a sufficient condition for the linear scattering of the Dirac equation.

\begin{lem}\label{s4: DL4.4}
Let $\psi$ be the global solution to the Dirac equation \eqref{s1: 1.1s}, where $\psi_0\in H^N(\mathbb{R}^2)$ for some $N\in\mathbb{N}$. Suppose that
\beq\label{s4: wd}
\int_0^{+\infty}\|F(\tau)\|_{H^N}{\rm{d}}\tau<\infty,
\eeq
then there exists some $\psi^+\in H^N(\mathbb{R}^2)$, satisfying
\beq\label{s4: psi-psi+}
\|\psi(t)-S(t)\psi^+\|_{H^N}\le C\int_t^{+\infty}\|F(\tau)\|_{H^N}{\rm{d}}\tau,\quad\mathrm{as}\ t\to+\infty,
\eeq
where $S(t):=e^{-t\gz^0\gz^a\partial_a}$ is the propagator for the linear Dirac equation, i.e, for any $\mathbb{C}^2$-valued function $f$,
\be
S(t)f=\mathcal{F}^{-1}e^{-it\gz^0\gz^a\xi_a}\mathcal{F}f(\xi).
\ee
Here $\mathcal{F}f$ is the Fourier transform of $f$ and $\mathcal{F}^{-1}$ is the inverse Fourier transform.

\begin{proof}
Using the Fourier transform, we can derive the following integral formula for \eqref{s1: 1.1s}
\beq\label{s4: ifpsi}
\psi(t,x)=S(t)\psi_0-i\int_0^tS(t-\tau)\gz^0F(\tau)\rm{d}\tau,
\eeq
where $S(t)$ verifies the following properties:
\begin{itemize}
\item[$i)$] $S(0)=I_2, S(t)S(\tau)=S(t+\tau)$;
\item[$ii)$] $\|\partial^I_xS(t)f\|_{L^2_x}=\|\partial^I_xf\|_{L^2_x}$, where $\partial^I_x:=\partial^{i_1}_{x_1}\partial^{i_2}_{x_2}$ for any multi-index $I=(i_1,i_2), i_1,i_2\in\mathbb{N}$, and hence $\|S(t)f\|_{H^N}=\|f\|_{H^N}$.
\end{itemize}
Let 
\be
\psi^+:=\psi_0-i\int_0^{+\infty}S(-\tau)\gz^0F(\tau){\rm{d}}\tau=\psi_0-i\lim_{T\to+\infty}\int_0^{T}S(-\tau)\gz^0F(\tau){\rm{d}}\tau.
\ee
We claim that $\psi^+$ is well-defined in $H^N(\mathbb{R}^2)$. Indeed, by the properties of $S(t)$ as above and the assumption \eqref{s4: wd}, we have
\ben
\l\|\int_{T_1}^{T_2}S(-\tau)\gz^0F(\tau){\rm{d}}\tau\r\|_{H^N}&\le&\int_{T_1}^{T_2}\|S(-\tau)\gz^0F(\tau)\|_{H^N}{\rm{d}}\tau\\
&\le&C\int_{T_1}^{T_2}\|F(\tau)\|_{H^N}{\rm{d}}\tau\to 0\quad\mathrm{as}\quad T_2>T_1\to+\infty.
\een
Hence the claim follows. The function $\varphi:=S(t)\psi^+$ solves the homogeneous Dirac equation $\mathcal{D}\varphi=0$ with $\varphi(0,\cdot)=\psi^+$. Furthermore, using the properties of $S(t)$ again, we obtain
\ben
\|\psi(t)-S(t)\psi^+\|_{H^N}&=&\l\|\lim_{T\to+\infty}\int_{t}^TS(t-\tau)\gz^0F(\tau)d\tau\r\|_{H^N}\\
&\le&\lim_{T\to+\infty}\int_{t}^T\|S(t-\tau)\gz^0F(\tau)\|_{H^N}\rm{d}\tau\\
&\le&C\int_t^{+\infty}\|F(\tau)\|_{H^N}\rm{d}\tau.
\een
\end{proof}
\end{lem}

\noindent$\textit{Proof\ of\ Theorem\ \ref{s1: thm2}.}$ Let $\psi$ be the global solution to \eqref{s1: 1.1s} obtained in Theorem \ref{s1: thm1} with $F(\psi)=(\psi^*\gz^0\psi)\psi$. By Lemma \ref{s4: DL4.4}, we only need to verify that for any $|I|\le N$,
\beq\label{s4: intFN}
\int_t^{+\infty}\|\partial^IF(\tau)\|_{L^2_x}{\rm{d}}\tau\le C\langle t\rangle^{-\f{1}{2}}\ln(2+t),
\eeq
and for any $|I|\le N-2$,
\beq\label{s4: intFN-2}
\int_t^{+\infty}\|\partial^IF(\tau)\|_{L^2_x}{\rm{d}}\tau\le C(t)\langle t\rangle^{-\f{1}{2}}\ln(2+t),\quad\quad\mathrm{with}\quad\quad\lim_{t\to+\infty}C(t)=0.
\eeq

Note that the proof of Theorem \ref{s1: thm1} in Sect. \ref{ss3.2} gives (see \eqref{s3: gzpsiinf}, \eqref{s3.2: gpsi-} and \eqref{s3.2: gpsi})
\ben
&&\sup_{t\ge 0}\l\{\sum_{|I|\le N}[E^D_{gst}(t,\hat{\Gamma}^I\psi)]^{\f{1}{2}}\r.\\
&+&\l.\sum_{|I|\le N-2}\l[\langle t+r\rangle^{\f{1}{2}}\langle t-r\rangle^{\f{1}{2}}|\Gamma^{I}\psi|+[\ln(2+t)]^{-1}\l(\langle t\rangle^{\f{3}{2}}|[\hat{\Gamma}^{I}\psi]_{-}|+\langle t-r\rangle\langle t\rangle^{\f{1}{2}}|\hat{\Gamma}^{I}\psi|\r)\r]\r\}\\
&\lesssim&C_1\ez.
\een

For any $|I|\le N$, using Lemma \ref{s2: DL4.2}, we have
\ben
|\partial^IF|&=&\l|\sum_{I_1+I_2+I_3=I}(\partial^{I_1}\psi)^*\gz^0(\partial^{I_2}\psi)(\partial^{I_3}\psi)\r|\\
&\lesssim&\sum_{|I_1|+|I_2|+|I_3|\le|I|}\l|[\partial^{I_1}\psi]_{-}^*\gz^0[\partial^{I_2}\psi]_{+}+[\partial^{I_1}\psi]_{+}^*\gz^0[\partial^{I_2}\psi]_{-}\r|\cdot|\partial^{I_3}\psi|\\
&\lesssim&\sum_{|I_1|+|I_2|+|I_3|\le|I|}|[\partial^{I_1}\psi]_{-}|\cdot|\partial^{I_2}\psi|\cdot|\partial^{I_3}\psi|,
\een
which implies
\beqn\label{s4: scFtau}
\|\partial^IF(\tau)\|_{L^2_x}&\lesssim&\sum_{\substack{|I_2|,|I_3|\le N-2,\\|I_1|\le N}}\|\langle r-\tau\rangle^{\f{1+\dz}{2}}\partial^{I_2}\psi\|_{L^\infty_x}\cdot\|\partial^{I_3}\psi\|_{L^\infty_x}\cdot\l\|\f{[\partial^{I_1}\psi]_{-}}{\langle r-\tau\rangle^{\f{1+\dz}{2}}}\r\|_{L^2_x}\nonumber\\
&+&\sum_{\substack{|I_1|,|I_2|\le N-2,\\|I_3|\le N}}\|[\partial^{I_1}\psi]_{-}\|_{L^\infty_x}\cdot\|\partial^{I_2}\psi\|_{L^\infty_x}\cdot\|\partial^{I_3}\psi\|_{L^2_x}\\
&\lesssim&(C_1\ez)^2\langle\tau\rangle^{-1}\ln(2+\tau)\sum_{|I'|\le N}\l\|\f{[\partial^{I'}\psi]_{-}}{\langle r-\tau\rangle^{\f{1+\dz}{2}}}\r\|_{L^2_x}+(C_1\ez)^3\langle\tau\rangle^{-2}\ln(2+\tau)\nonumber.
\eeqn
It follows that
\ben
\int_t^{+\infty}\|\partial^IF(\tau)\|_{L^2_x}\rm{d}\tau&\lesssim&(C_1\ez)^2\sum_{|I'|\le N}\l(\int_t^{+\infty}\langle\tau\rangle^{-2}\ln^2(2+\tau) \rm{d}\tau\r)^{\f{1}{2}}\l(\int_t^{+\infty}\l\|\f{[\partial^{I'}\psi]_{-}}{\langle r-\tau\rangle^{\f{1+\dz}{2}}}\r\|_{L^2_x}^2\rm{d}\tau\r)^{\f{1}{2}}\\
&+&(C_1\ez)^3\int_t^{+\infty}\langle\tau\rangle^{-2}\ln(2+\tau) \rm{d}\tau\\
&\lesssim&(C_1\ez)^3\langle t\rangle^{-\f{1}{2}}\ln(2+t),
\een
which implies \eqref{s4: intFN}.

It remains to prove \eqref{s4: intFN-2}. For $|I|\le N-2$, we can omit the second sum in \eqref{s4: scFtau} and derive
\ben
\int_t^{+\infty}\|\partial^IF(\tau)\|_{L^2_x}{\rm{d}}\tau&\lesssim&(C_1\ez)^2\sum_{|I'|\le N}\l(\int_t^{+\infty}\langle\tau\rangle^{-2}\ln^2(2+\tau) \rm{d}\tau\r)^{\f{1}{2}}\l(\int_t^{+\infty}\l\|\f{[\partial^{I'}\psi]_{-}}{\langle r-\tau\rangle^{\f{1+\dz}{2}}}\r\|_{L^2_x}^2\rm{d}\tau\r)^{\f{1}{2}}\\
&\lesssim& C(t)(C_1\ez)^2\langle t\rangle^{-\f{1}{2}}\ln(2+t),
\een
where
\be
C(t):=\sum_{|I'|\le N}\l(\int_t^{+\infty}\l\|\f{[\partial^{I'}\psi]_{-}}{\langle r-\tau\rangle^{\f{1+\dz}{2}}}\r\|_{L^2_x}^2{\rm{d}}\tau\r)^{\f{1}{2}}\to 0\quad\mathrm{as}\ t\to+\infty.
\ee

\begin{appendices}

\section{}\label{sB}
In this section we show that the assumption on the initial data in Theorem \ref{s1: thm1} can be relaxed.

Let $u$ be the solution to 
\beq\label{sB: uh}
\l\{\begin{array}{rcl}
-\Box u(t,x)&=&f(t,x),\\
(u,\partial_tu)|_{t=0}&=&(u_0,u_1).
\end{array}\r.
\eeq
Denote
\ben
\mathcal{E}(t,u):&=&\|\partial u(t,x)\|_{L^2_x}^2,\quad\quad|\partial u|^2=\sum_{\az}(\partial_\az u)^2,\\
\mathcal{E}_{con}(t,u):&=&\sum_{a=1}^2\|(|L_0 u+u|+|\Omega_{12} u|+|L_a u|)(t,x)\|_{L^2_x}^2.
\een
\begin{lem}(See \cite{A, LZ}.)\label{sB: hw}
Let $u$ be the solution to \eqref{sB: uh}. Then the following estimates hold:
\begin{itemize}
\item[$i)$] (Standard energy estimate)
$$\mathcal{E}(t,u)^{\f{1}{2}}\lesssim \mathcal{E}(0,u)^{\f{1}{2}}+\int_0^t\|f(\tau,x)\|_{L^2_x}\rm{d}\tau;$$
\item[$ii)$] (Conformal energy estimate)
$$\mathcal{E}_{con}(t,u)^{\f{1}{2}}\lesssim \mathcal{E}_{con}(0,u)^{\f{1}{2}}+\int_0^t\|\langle \tau+|x|\rangle f(\tau,x)\|_{L^2_x}\rm{d}\tau;$$
\item[$iii)$] ($L^2$ estimate) 
$$\|u(t,x)\|_{L^2_x}\lesssim\|u_0\|_{L^2_x}+\ln^{\f{1}{2}}(2+t)\l[\|u_1\|_{L^1_x\cap L^2_x}+\int_0^t\|f(\tau,x)\|_{L^1_x\cap L^2_x}\rm{d}\tau\r].$$
\end{itemize}
\end{lem} 

\begin{lem}\label{sB: gu}
Let $u$ be the solution to \eqref{sB: uh} with $f\equiv 0$. Then for any integer $N\ge 3$, we have
\beqn\label{sB: gammaIu}
\sum_{|I|\le N-2}\|\Gamma^Iu\|_{L^\infty_x}&\lesssim&\langle t\rangle^{-\f{1}{2}}\Bigg\{\|u_0\|_{L^2_x}+\ln^{\f{1}{2}}(2+t)\|u_1\|_{L^1_x\cap L^2_x}\nonumber\\
&+&\sum_{k\le N}\|\langle|x|\rangle^{k}\nabla^ku_0\|_{L^2_x}+\sum_{k\le N-1}\|\langle|x|\rangle^{k+1}\nabla^ku_1\|_{L^2_x}\Bigg\}.
\eeqn
\begin{proof}
By Lemma \ref{sB: hw}, we have
\beq\label{sB: uL2}
\|u(t,x)\|_{L^2_x}\lesssim\|u_0\|_{L^2_x}+\ln^{\f{1}{2}}(2+t)\|u_1\|_{L^1_x\cap L^2_x}
\eeq
and
\be
\mathcal{E}(t,u)^{\f{1}{2}}+\mathcal{E}_{con}(t,u)^{\f{1}{2}}\lesssim \mathcal{E}(0,u)^{\f{1}{2}}+\mathcal{E}_{con}(0,u)^{\f{1}{2}},
\ee
which implies
\beqn\label{sB: gamma1}
\sum_{|I|\le 1}\|\Gamma^Iu\|_{L^2_x}&\lesssim&\|u\|_{L^2_x}+\sum_{|I|\le 1}\|\Gamma^Iu(0)\|_{L^2_x}\nonumber\\
&\lesssim&\|u\|_{L^2_x}+\sum_{k\le 1}\|\langle|x|\rangle^{k}\nabla^ku_0\|_{L^2_x}+\sum_{k=0}\|\langle|x|\rangle^{k+1}\nabla^ku_1\|_{L^2_x}.
\eeqn
Acting the vector fields $\Gamma_k,k=1,\cdots,7$ on both sides of \eqref{sB: uh} and apply the standard energy and conformal energy estimates in Lemma \ref{sB: hw}, we have
\be
\mathcal{E}(t,\Gamma_ku)^{\f{1}{2}}+\mathcal{E}_{con}(t,\Gamma_ku)^{\f{1}{2}}\lesssim \mathcal{E}(0,\Gamma_ku)^{\f{1}{2}}+\mathcal{E}_{con}(0,\Gamma_ku)^{\f{1}{2}},
\ee
which combined with \eqref{sB: gamma1} gives
\ben\label{sB: gamma2}
\sum_{|I|\le 2}\|\Gamma^Iu\|_{L^2_x}&\lesssim&\|u\|_{L^2_x}+\sum_{|I|\le 2}\|\Gamma^Iu(0)\|_{L^2_x}\nonumber\\
&\lesssim&\|u\|_{L^2_x}+\sum_{k\le 2}\|\langle|x|\rangle^{k}\nabla^ku_0\|_{L^2_x}+\sum_{k\le 1}\|\langle|x|\rangle^{k+1}\nabla^ku_1\|_{L^2_x}.
\een
By induction, for any $N\in\mathbb{N}$ with $N\ge 1$, we have
\beqn\label{sB: gammaN}
\sum_{|I|\le N}\|\Gamma^Iu\|_{L^2_x}&\lesssim&\|u\|_{L^2_x}+\sum_{|I|\le N}\|\Gamma^Iu(0)\|_{L^2_x}\nonumber\\
&\lesssim&\|u\|_{L^2_x}+\sum_{k\le N}\|\langle|x|\rangle^{k}\nabla^ku_0\|_{L^2_x}+\sum_{k\le N-1}\|\langle|x|\rangle^{k+1}\nabla^ku_1\|_{L^2_x}\nonumber\\
&\lesssim&\|u_0\|_{L^2_x}+\ln^{\f{1}{2}}(2+t)\|u_1\|_{L^1_x\cap L^2_x}+\sum_{k\le N}\|\langle|x|\rangle^{k}\nabla^ku_0\|_{L^2_x}+\sum_{k\le N-1}\|\langle|x|\rangle^{k+1}\nabla^ku_1\|_{L^2_x},\nonumber\\
\eeqn
where we use \eqref{sB: uL2} in the last inequality. By Lemma \ref{s2: K-S}, for any $|I|\le N-2$, we have
\ben
\|\Gamma^Iu\|_{L^\infty_x}
&\lesssim&\langle t\rangle^{-\f{1}{2}}\sum_{|I|\le N}\|\Gamma^Iu\|_{L^2_x}\nonumber\\
&\lesssim&\langle t\rangle^{-\f{1}{2}}\l\{\|u_0\|_{L^2_x}+\ln^{\f{1}{2}}(2+t)\|u_1\|_{L^1_x\cap L^2_x}+\sum_{k\le N}\|\langle|x|\rangle^{k}\nabla^ku_0\|_{L^2_x}+\sum_{k\le N-1}\|\langle|x|\rangle^{k+1}\nabla^ku_1\|_{L^2_x}\r\},
\een
hence \eqref{sB: gammaIu} follows.
\end{proof}
\end{lem}

We recall \eqref{s3: DL5.3'} and the estimate \eqref{s3: gPsiB} in $\mathbf{Step\ 2}$ in the proof of Proposition \ref{s3: prop 1}. We decompose $\Psi$ into the homogeneous wave component $\Psi_1$ and the inhomogeneous wave component $\Psi_2$, and use \eqref{sB: gammaIu} for $\Gamma^I\Psi_1$ and Lemma \ref{s2: DMY2.78} for $\Gamma^I\Psi_2$. For any $|I|\le N-2$, we obtain
\beqn\label{sB: gPsi}
\|\Gamma^I\Psi(t)\|_{L^\infty_x}&\lesssim&\langle t\rangle^{-\f{1}{2}}\Bigg\{\ln^{\f{1}{2}}(2+t)\|\psi_0\|_{L^1_x\cap L^2_x}+\sum_{k\le N-1}\|\langle|x|\rangle^{k+1}\nabla^k\psi_0\|_{L^2_x}\nonumber\\
&+&\|\Gamma^I\Psi_2(0,\cdot)\|_{W^{2,1}}+\|\partial_t\Gamma^I\Psi_2(0,\cdot)\|_{W^{1,1}}+\sum_{|I'|\le|I|+1}\int_0^t\langle\tau\rangle^{-\f{1}{2}}\|\Gamma^{I'}F(\tau)\|_{L^1_x}\rm{d}\tau\Bigg\}.\nonumber\\
\eeqn
Then following the proof of Proposition \ref{s3: prop 1} and choosing $N$ larger ($N\ge 5$), we can obtain Theorem \ref{s1: thm1} with relaxed condition on the smallness of the initial data, i.e.,
\be
\|\psi_0\|_{L^1_x}+\sum_{k\le N}\|\langle|x|\rangle^{k+1}\nabla^k\psi_0\|_{L^2_x}\le\ez.
\ee

\begin{rem}
In \eqref{sB: gPsi}, we note that initial data $(\Gamma^I\Psi_2,\partial_t\Gamma^I\Psi_2)|_{t=0}\neq (0,0)$. However, since the nonlinearity $F$ is cubic (see \eqref{s3: DL5.3'}), using the H\"older inequality, we can obtain 
\be
\sum_{|I|\le N-2}\l(\|\Gamma^I\Psi_2(0,\cdot)\|_{W^{2,1}}+\|\partial_t\Gamma^I\Psi_2(0,\cdot)\|_{W^{1,1}}\r)\lesssim\sum_{k\le N}\|\langle|x|\rangle^{k+1}\nabla^k\psi_0\|_{L^2_x}.
\ee
\end{rem}

\end{appendices}

\bigskip

\section*{Data Availability Statement}

Data sharing is not applicable to this article as no datasets were generated or analysed during the current study.

Qian Zhang

School of Mathematical Sciences, Beijing Normal University

Laboratory of Mathematics and Complex Systems, Ministry of Education Beijing 100875, China

Email: qianzhang@bnu.edu.cn


\begin{thebibliography}{99}

\bibitem{Al1} \label{Al1} Alinhac, S.:
\textit{The null condition for quasilinear wave equations in two space dimensions I,}
Invent. Math. $\mathbf{145}$ (3) (2001), 597-618.

\bibitem{A} \label{A} Alinhac, S.:
\textit{Hyperbolic Partial Differential Equations.}
Springer-Verlag, New York, 2009.

\bibitem{Ba} \label{Ba} Bachelot, A.:
\textit{Probl\`eme de Cauchy global pour des syst\`emes de Dirac-Klein-gordon.}
Ann. Inst. H. Poincar\'e Phys. Th\'eor. $\mathbf{48}$ (4) (1988), 387-422.

\bibitem{BH} \label{BH} Bejenaru, I., Herr, S.: 
\textit{The cubic Dirac equation: small initial data in $H^1(\mathbb{R}^3)$.} 
Comm. Math. Phys. $\mathbf{335}$ (1), (2015), 43-82.

\bibitem{BH2} \label{BH2} Bejenaru, I., Herr, S.: 
\textit{The cubic Dirac equation: small initial data in $H^{\f{1}{2}}(\mathbb{R}^2)$.} 
Comm. Math. Phys., $\mathbf{343}$ (2) (2016), 515-562. 


\bibitem{Bo} Bournaveas, N.:
\textit{Low regularity solutions of the Dirac-Klein-Gordon equations in two space dimensions.} 
Comm. Partial Differential Equations $\mathbf{26}$ (7-8) (2001), 1345-1366.

\bibitem{BC} Bournaveas, N., Candy, T.:
\textit{Global well-posedness for the massless cubic Dirac equation.} 
Int. Math. Res. Not. IMRN 2016, no. 22, 6735-6828.

\bibitem{D21} \label{D21} Dong, S.,
\textit{Global solution to the wave and Klein-Gordon system under null condition in dimension two.} 
J. Funct. Anal. $\mathbf{281}$ (11) (2021), Paper No. 109232, 29 pp.

\bibitem{DL} \label{DL} Dong, S., Li, K.:
\textit{Global solution to the cubic Dirac equation in two space dimensions.} 
arXiv:2111.04048

\bibitem{DLMY} \label{DLMY} Dong, S. Li, K., Ma, Y., Yuan, X.:
\textit{Global behavior of small data solutions for the $2D$ Dirac-Klein-Gordon Equations.} 
arXiv:2205.12000

\bibitem{DLW} \label{DLW} Dong, S., LeFloch, P.G., Wyatt, Z.:
\textit{Global evolution of the U(1) Higgs Boson: nonlinear stability and uniform energy bounds.}
Ann. Henri Poincar\'e $\mathbf{22}$ (3) (2021), 677-713.

\bibitem{DW} \label{DW} Dong, S., Wyatt, Z.:
\textit{Hidden structure and sharp asymptotics for the Dirac-Klein-Gordon system in two space dimensions.}
arXiv:2105.13780


\bibitem{Ev} \label{Ev} Escobedo, M., Vega, L.:
\textit{A semilinear Dirac equation in $H^s(\mathbb{R}^3)$ for $s>1$.}
SIAM J. Math. Anal. $\mathbf{28}$ (2) (1997), 338-362.

\bibitem{Ge} \label{Ge} Georgiev, V.:
\textit{Decay estimates for the Klein-Gordon equation.}
Comm. Partial Differential Equations $\mathbf{17}$ (7-8) (1992), 1111-1139.

\bibitem{Ho} \label{Ho} H\"ormander, L.:
\textit{Lectures on nonlinear hyperbolic differential equations.}
Math\'ematiques \& Applications (Berlin) [Mathematics \& Applications], vol. 26, Springer-Verlag, Berlin, 1997, viii+289 pp.

\bibitem{J81}\label{J81} John F.:
\textit{Blow up of solutions for quasi-linear wave equations in three space dimensions,} 
Commun. Pure Appl. Math. $\mathbf{34}$ (1) (1981), 29-51.


\bibitem{K85} \label{K85} Klainerman, S.:  
\textit{Uniform decay estimates and the Lorentz invariance of the classical wave equation,} 
Comm. Pure Appl. Math. $\mathbf{38}$ (3) (1985), 321-332.

\bibitem{K86} \label{K86} Klainerman, S.:
\textit{The null condition and global existence to nonlinear wave equations.}
Nonlinear systems of partial differential equations in applied mathematics, Part 1 (Santa Fe, N.M., 1984), 293-326, Lectures in Appl. Math., vol. 23, Amer. Math. Soc., Providence, RI, 1986.

\bibitem{LZ} \label{LZ} Li, T., Zhou, Y.:
\textit{Nonlinear wave equations.}
Vol. 2. Translated from the Chinese by Yachun Li. Series in Contemporary Mathematics, 2. Shanghai Science and Technical Publishers, Shanghai; Springer-Verlag, Berlin, 2017. xiv+391 pp.

\bibitem{MNN} \label{MNN} Machihara, S., Nakamura, M., Nakanishi, K., Ozawa, T.:
\textit{Endpoint Strichartz estimates and global solutions for the nonlinear Dirac equation.}
J. Funct. Anal. $\mathbf{219}$ (1) (2005), 1-20.

\bibitem{P} \label{P} Pecher, H.: 
\textit{Local well-posedness for the nonlinear Dirac equation in two space dimensions.} 
 Commun. Pure Appl. Anal. $\mathbf{13}$ (2), (2014), 673-685.

\bibitem{S} \label{S} Sogge, C.D.:
\textit{Lectures on non-linear wave equations.} 
Second edition. International Press, Boston, MA, 2008. x+205 pp.

\bibitem{So} \label{So} Soler, M.:
\textit{Classical, stable, nonlinear spinor field with positive rest energy.}
Phys. Rev. D $\mathbf{1}$ (1970), 2766-2769.

\bibitem{T} \label{T} Thirring, W.E.:
\textit{A soluble relativistic field theory.}
Ann. Physics $\mathbf{3}$ (1) (1958), 91-112.

\bibitem{Tz} \label{Tz} Tzvetkov, N.: 
\textit{Existence of global solutions to nonlinear massless Dirac system and wave equation with small data.} 
Tsukuba J. Math. $\mathbf{22}$ (1) (1998), 193-211.

\end{thebibliography}
\end{document}